\documentclass[11pt]{article} 

\usepackage{graphicx}
\usepackage{amsfonts,amsmath,graphicx,color,epsfig}
\usepackage{todonotes}
\usepackage{changes}
\usepackage[font=small,labelfont=bf]{caption}
\usepackage{subcaption}
\usepackage{amsthm}
\usepackage{amssymb}
\usepackage{latexsym}
\usepackage{color}
\usepackage{mathrsfs}
\usepackage{natbib}
\usepackage{hyperref}
\usepackage[toc,page]{appendix}
\usepackage[left=2.5cm,top=4cm,right=2.5cm,bottom=3cm]{geometry}
\usepackage{float}
\usepackage{subfloat}

\usepackage{colortbl,xcolor}

\newcommand{\dW}{{\partial W}}

\newcommand{\dd}{\, \text{d}}

\newcommand{\bP}{{\mathbb{P}}}

\newcommand{\bR}{{\mathbb{R}}}
\newcommand{\bE}{{\mathbb{E}}}

\newcommand{\lck}{\left\lbrack}
\newcommand{\rck}{\right\rbrack}
\newcommand{\lce}{\left\lbrace}
\newcommand{\rce}{\right\rbrace}

\newcommand{\un}{\mathbb{I}}

\newtheorem{proposition}{Proposition}

\definecolor{Gray}{gray}{0.85}

\usepackage{setspace}
\begin{document}

\begin{center}
\large{\textbf{Note on the approximation of the conditional intensity }}

\large{\textbf{of non-stationary cluster point processes}}

\bigskip

Edith \textsc{Gabriel}\footnote{\tt{edith.gabriel@inrae.fr}} and Jo{\"e}l \textsc{Chad{\oe}uf}

\bigskip

 Biostatistics and Spatial Processes Unit, INRAE, F-84911 Avignon, France

\bigskip

\today

\vspace{1cm}

\begin{minipage}{12.5cm}
  {\small
\textbf{Abstract}

In this note we consider non-stationary cluster point processes and we derive their conditional intensity, i.e. the intensity of the process given the locations of one or more events of the process. We then provide some approximations of the conditional intensity.

\medskip

\textsc{Conditional intensity}; \textsc{Neyman-Scott process}; \textsc{Point process}
}
\end{minipage}
\end{center}

\vspace{.5cm}

\section{Introduction}

The problem of conditioning is an old problem in the point process theory. The conditional distribution of a spatial point process, say $\Phi$, given a realisation of some events of $\Phi$ was introduced by \cite{palm1943} for stationary point processes on the real line and recently summed up in \cite{coeurjolly2017} in the general case.

We consider locally finite point processes $\Phi$ defined in a compact set $S \subset \bR^2$ and specified by a density~$f$.
We assume that we have observed $\Phi$ in $W \subseteq S$ and we denote $\Phi_W$ its restriction to the set~$W$. A natural way to predict $\Phi_{S \backslash W}$ is to consider the conditional distribution of $\Phi_{S \backslash W}$ given $\Phi_W$, which is can be expressed in terms of the conditional density (cf \cite{coeurjolly2017}). Unfortunately, the density of $\Phi$ restricted to $W$ is tractable for few usual processes, as Poisson, Gibbs and determinantal processes, but not for Cox and cluster point processes.

In this note, we derive the conditional intensity of non-stationary cluster point processes and provide some approximations for practical applications.
Note that  \cite{gabriel2017} and  \cite{gabriel2021} define a ``model-free'' predictor of the conditional intensity for stationary and non-stationary processes, in sense that it is only related to the first and second-order moments of the point process.

\section{Conditional intensity of non-stationary cluster point processes}

Cluster point processes, developed by \cite{neymanscott1958}, are formed by a simple procedure, with a homogeneous Poisson process $\Psi$ with intensity $\kappa$ generating \textit{parent} points  at a first step and a random pattern of \textit{offspring} points around each parent point at a second step.
The number of offspring points has a Poisson distribution with mean $\mu$ and the offspring points are independently and identically distributed with a bounded support kernel $k$ depending on the distance from offspring to parent.
The cluster point process $\tilde \Phi$ is the set of offspring points, regardless their parentage. This process is stationary with intensity $\tilde \lambda = \kappa \mu$.

\medskip
We focus on the $p(x)-$thinned process of  $\tilde \Phi$, where $p(x)$ is a deterministic function on $\bR^2$ with $0 \leq p(x) \leq 1$.  If the point $x$
belongs to $\tilde \Phi$, it is deleted with probability $1-p(x)$ and again its deletion is independent of locations and possible deletions of any other points.
Let $\Phi$ be the $p(x)-$thinned process.
The process $\Phi$ is second-order intensity reweighted stationary with $\lambda (x) = \kappa \mu p(x)$ (see \cite{chiu2013}).

\medskip

\begin{figure}[h]
\centering
  \includegraphics[width=0.25\linewidth]{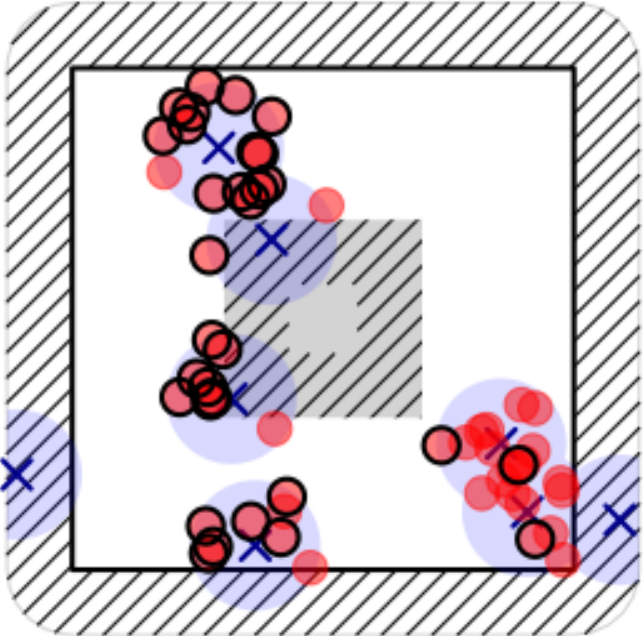}
\caption{Generating procedure of the $p(x)$-thinned cluster process. Parent points ($\Psi=$blue crosses) are generated in the union of the observation window $W$ (white area) and $\dW$ (hatched area). Offspring ($\tilde \Phi=$red dots) are generated within the circle around the parent points (blue area). The final process ($\Phi=$black dots) is obtained by $p(x)$-thinning. The window of interest $S$ is delineated in black. The prediction window is represented by the grey shaded central square.}
\label{fig:clusterprocess}
\end{figure}

Here we want to know the conditional intensity of $\Phi_{S \backslash W}$ given $\Phi_W$.
We denote $\dW$ the border of the observation window $W$, with width defined by the range of the dispersion kernel $k$, say $r$. In other words, $\dW = W_{\oplus r} \backslash W$. Figure~\ref{fig:clusterprocess} illustrates the different steps of the generating procedure.

\begin{proposition}
\label{prop:intcond}
The conditional intensity of a $p(x)$-thinned cluster process observed in $W \subseteq S$ is:
\begin{multline}
\lambda (x_o |  \Phi_W) = \int \left\lbrack \sum_{y \in W \cup \dW} \mu p(x_o) k(y-x_o)
 \right. \\
\left. + \mu \kappa \int_{b(x_o,r) \backslash (W \cup \dW)} p(x_o) k(y-x_o) \dd y \right\rbrack
\dd \bP[ \Psi_{W \cup \dW} |  \Phi_W],  \label{eq:intcond}
\end{multline}
where $b(x_o,r)$ denotes the disc of centre $x_o$ and radius $r$.
\end{proposition}
\noindent
The proof is rather straightforward as $(i)$ for a cluster process we know the conditional intensity given the realisation of parent points, $(ii)$ the parent process is Poisson and $(iii)$ knowing the offspring points in $W$ is not informative on parent points in $b(x_o,r) \backslash (W \cup \dW)$.

\medskip
\noindent
The conditional intensity~(\ref{eq:intcond}) depends on the distribution of the parent process in $W \cup \dW$ given the offspring in $W$.
 By Campbell's theorem this equation can be rewritten
\begin{multline}
\lambda (x_o |  \Phi_W) = \int_{W \cup \dW} \mu p(x_o) k(y-x_o)  \rho(y |\Phi_W) \dd y \\
 + \mu \kappa \int_{b(x_o,r) \backslash (W \cup \dW)} p(x_o) k(y-x_o) \dd y,
 \label{eq:intcond2}
\end{multline}
 where $\rho(y |\Phi_W)$ is the intensity of parent points in $W \cup  \dW$ given the offspring points in $W$.
 \cite{baudin1983} derived the following formula for such a conditional intensity of parent points for Neyman-Scott processes:
 \vspace{-4mm}
 \begin{multline}
\rho(y | \Phi_W) = \kappa G \left(1-F(W_{y}) \right)
+ \sum_{j=1}^{2^n-1} \sum_{b \in {\cal B}} b(a_j) \prod_{i=1}^{2^n-1} S(\Phi,W,a_i)^{b(a_i)} \\
\times \left\lbrack \dfrac{\kappa G^{|a_j|}\left(1-F(W_{y}) \right) \prod_{\ell =1}^n k(y_\ell - x_o)^{a_{j\ell}}}{S(\Phi,W,a_j)} \right\rbrack^{b(a_j)}  \left\lbrack \sum_{b \in {\cal B}} \prod_{i=1}^{2^n-1} S(\Phi,W,a_i)^{b(a)} \right\rbrack,
\label{eq:condparent}
\end{multline}
where
\begin{itemize}
\itemsep0em
    \item  [{\tiny $\bullet$}] $G$ is the probability generating function of the number of points in a cluster,
    \item [{\tiny $\bullet$}] $F(\dd x)$ is the probability distribution function of offspring points, with density $k$,
   \item [{\tiny $\bullet$}] $\lce x_1, \dots, x_n \rce = \Phi_W$,
    \item [{\tiny $\bullet$}] $W_{y} = - y + W$,
    \item [{\tiny $\bullet$}] $a_1=(0,\dots,0,1)$, $a_2=(0,\dots,0,1,0)$, \dots, $a_{2^n-1}=(1,\dots,1)$: vectors of length $n$,
    \item [{\tiny $\bullet$}] ${\cal B}$ is the set of all functions $b$ : $\lce a_1, \dots, a_{2^n-1}\rce \to \lce 0,1 \rce$ such that $\sum_{i=1}^{2^n-1} b(a_i) a_i = (1,\dots,1)$, $a_i=(a_{i1},\dots,a_{in})$, $|a_i| = a_{i1} + \dots + a_{in}$,
   \item [{\tiny $\bullet$}] $S(\Phi,W,a_i) = \kappa \int G^{|a_i|}(1 - F(W_y)) \prod_{\ell=1}^n k(x_l -y) \dd y.$
\end{itemize}
 However, this conditional intensity is based on combinations and is just too complicated in practice. \cite{vanLieshout1995} interpreted the problem of identifying parent points as a statistical estimation problem with a Bayesian inference based on MCMC methods; see also \cite{lawson2002} for similar approaches.

\medskip

Here we propose to approximate it as follows,
\begin{equation}
\rho (y |  \Phi_W) = \frac{c(y)}{\mu p(y)} \sum_{x \in \Phi_W} k(x-y)
 +  \kappa \exp \left( - \mu \int_{W} p(z) k(y-z) \dd z \right)
 \label{eq:intcondpar}
\end{equation}
where $c(y)$ ensures that $\bE \lck \rho (y |  \Phi_W)  \rck = \kappa$.
The first term of Equation~(\ref{eq:intcondpar}) is related to observed offspring points in $W$ and the second term to offspring points that have not been observed (due to thinning or their proximity to the boundary of $W$).

\section{Validation procedure}

In this section, we aim at quantifying the difference between the distribution of parent points given the offspring points provided in~(\ref{eq:intcondpar}) and the observed distribution of the conditional parent process. We thus use complementary statistics to test the interactions $(i)$ between parent points, $(ii)$ between parent and offspring points and $(iii)$ between parent points and the inner (or outer) boundary of $W$. Let $\Phi(A)$ (resp. $\Phi(B)$) be the number of points of $\Phi$ in a Borel set $A$ (resp. of $\Psi$ in $B$). We denote $\Psi_S$ the observed parent points in $S$ and $\nu(B)$ the Lebesgue measure of $B$. Then,

\noindent
\textit{(i) Interaction statistic between parent points}

Let $\widehat H(d)$ be the empirical cumulative distribution function between observed parent points $\Psi_S$ and $H(d)$ the theoretical one:
\begin{eqnarray*}
\widehat H(d) & = & \dfrac{1}{\Psi(S)} \sum_{y_1,y_2 \in \Psi_S}^{\neq} \un_{\lce \| y_1-y_2 \| \leq d \rce},  \\
H(d) & = & \dfrac{1}{\int_S \rho(y | \Phi_W) \dd y} \int_S \int_{b(z,d) \cap S} \rho(y | \Phi_W) \rho(z | \Phi_W) \dd y \dd z.
\end{eqnarray*}

\noindent
 \textit{(ii) Interaction statistic between parent and offspring points}

Let $\widehat E(d)$ be the empirical cumulative distribution function between observed parent points $\Psi_S$ and observed offspring points $\Phi_W$ and $E(d)$ the theoretical one:
\begin{eqnarray*}
\widehat E(d) & = & \dfrac{1}{\Phi(W)} \sum_{x \in \Phi_W} \sum_{y \in \Psi_S} \un_{\lce \| x-y \| \leq d \rce},\\
E(d) & = & \dfrac{1}{\Phi(W)}  \sum_{x \in \Phi_W} \int_{b(z,d) \cap S} \rho(z | \Phi_W) \dd z.
\end{eqnarray*}

\noindent
\textit{(iii) Interaction statistic between parent points and the boundary $b_W$ of $W$}

Let $\widehat B(d)$ be the empirical cumulative distribution function between observed parent points $\Psi_S$ and the inner boundary ($B^{inner}$)  or outer boundary ($B^{outer}$) of $W$, and $B(d)$ the theoretical one:
\begin{eqnarray*}
\widehat B(d) & = & \dfrac{1}{\nu(b_W)} \sum_{\ell \in b_W} \sum_{y \in \Psi_S} \un_{\lce \| \ell-y \| \leq d \rce} \\
& = &  \dfrac{1}{\nu(b_W)} \sum_{y \in \Psi_S} \nu(b_W \cap b(y,d)), \\
B (d) & = & \dfrac{1}{\nu(b_W)} \int_{b_W} \int_{b(\ell,d) \cap S} \rho(z | \Phi_W) \dd z \dd \ell.
\end{eqnarray*}

\medskip
We illustrate the results for a thinned Mat\'ern cluster process $\Phi$. For this process, $k$ is the uniform distribution on the disc of radius $r$ and the conditional intensity is
\begin{equation}
     \lambda (x_o |  \Phi_W) = \frac{\mu p(x_o)}{\pi R^2} \int_{b(x_o,r) \cap (W \cup \dW)} \rho (y | \Phi_W) \dd y
   + \kappa \mu p(x_o) \nu \left( b(x_o,r) \backslash (W \cup \dW) \right),
\end{equation}
with
\begin{equation*}
  \rho (y |  \Phi_{W}) =  \frac{1}{\mu p(y) \pi R^2} \sum_{x \in \Phi_{W}} \un_{b(x,r)} (y)
  + \kappa \exp \left( - \frac{\mu}{\pi r^2} \int_{b(y,r) \cap W} p(z) \dd z \right).
\end{equation*}
The non-stationary Mat{\'e}rn Cluster process  $\Phi$ depends on four parameters: the thinning probability $p(x)$, the intensity of parents $\kappa$, the mean number of points per parent $\mu$ and the radius of dispersion of the offspring around the parent points $r$.
Here we fix $\kappa=50$ and $\mu=40$ and we consider:

\noindent
- two thinning probabilities:
  $p_1(x)=p_1(x_1,x_2)=\alpha_1 \un_{\lce x_1 \leq v \rce} + \alpha_2 \un_{\lce x_1 > v \rce}$, setting $\alpha_1 = 0.8$, $\alpha_2=0.2$ and $v=0.5$,
  and $p_2(x) = p_2(x_1,x_2) =  1-x_1$.

\noindent
- the unit square as study region $S$. The observation window is $W = S \backslash W_h$, where $W_{h} = \lbrack 0.35, 0.65 \rbrack^2$ when using $p_1(x)$ and
        $W_{h} = \lbrack 0.05, 0.95 \rbrack \times \lbrack 0.36, 0.64 \rbrack  $ when using $p_2(x)$.

\noindent
- $r \in \lbrace 0.05, 0.09, 0.13\rbrace$.

\begin{figure}[h!]
\centering
\begin{subfigure}{.45\textwidth}
  \centering
  \includegraphics[width=\linewidth]{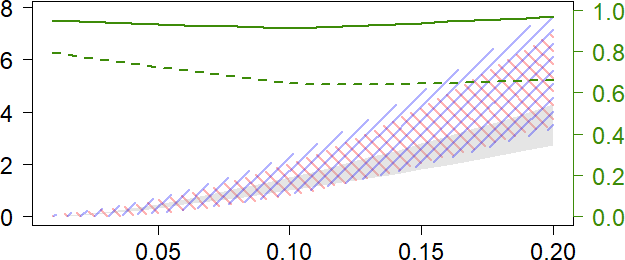}
  \caption{Interactions between parent points.}
\end{subfigure}%

\begin{subfigure}{.45\textwidth}
  \centering
  \includegraphics[width=\linewidth]{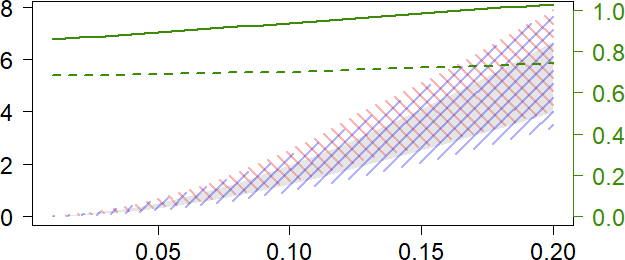}
    \caption{Interactions between parent and offspring points.}
\end{subfigure}


\begin{subfigure}{.45\textwidth}
  \centering
  \includegraphics[width=\linewidth]{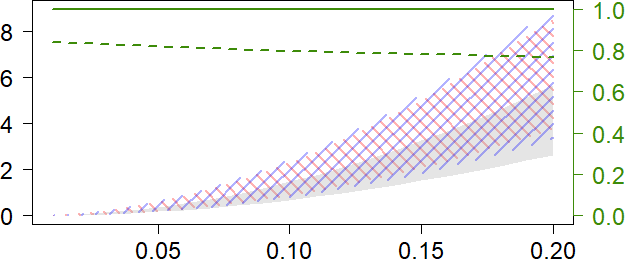}
    \caption{Interactions between parent points and the inner boundary of $W$.}
\end{subfigure}

\begin{subfigure}{.45\textwidth}
  \centering
  \includegraphics[width=\linewidth]{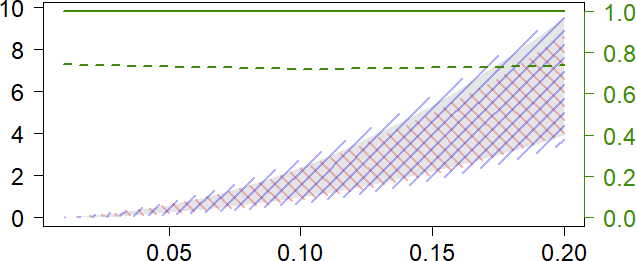}
    \caption{Interactions between parent points and the outer boundary of $W$.}
\end{subfigure}
\caption{95\% envelopes of the cumulative distribution functions calculated for $N$ realisations of the thinned Mat{\'e}rn cluster process (empirical is in red hatching and theoretical is in grey) and for simulated parent points (blue hatching). The curves represent the cover rates $\tau_1$ (solid) and $\tau_2$ (dashed).
}
\label{fig:inter_009}
\end{figure}

\noindent
For each pair of parameters $(p(x),r)$ we simulate $N=250$ realisations of the non-stationary Mat{\'e}rn Cluster process and compute all the previous interaction statistics.
Because all the theoretical statistics only provide a trend, for each of the $N$ simulations, we generate $n=100$ simulations of parent points from a Poisson process with intensity $\rho (y| \Phi_W)$ and compute the related empirical cumulative distribution functions, that we denote by $\widehat H_{sim}(d)$, $\widehat E_{sim}(d)$ and $\widehat B_{sim}(d)$.
Figure~\ref{fig:inter_009} illustrate the 95\% envelopes of the empirical cumulative distribution functions computed from the $N$ observed parent points (red hatching) and from the $N \times n$ simulated parent points (blue hatching). The grey envelopes correspond to the theoretical trend. In this figure $p(x) = p_1(x)$ and $r=0.09$.

\medskip

\noindent
Results for all pairs of parameters are very similar. All overlapping envelopes indicate that the empirical cumulative distribution functions between parent points and other parent points / offspring points and the boundary of the observation window $W$ are similar for the observed parent points $\Psi_S$ and for simulated parent points, which further correspond to the theoretical distribution. This is true at any distances and shows that the
main characteristics of the approximated distribution of parent points in $S$ given the offspring points in $W$ provided in~(\ref{eq:intcondpar}) include those of the original distribution of parent points in $S$.

\medskip

\noindent
For each type of interaction, we computed the global coverage rates between the envelopes obtained from the observed distribution of parent points and the enveloped obtained from the approximated distribution of parent points. E.g., denoting by ${\cal E}$ the envelopes, the coverage rates for the interaction statistic between parent points are
 $$\tau_1 (H) = \nu \left( {\cal E}(\widehat H) \cap {\cal E}(\widehat H_{sim}) \right) / \nu \left( {\cal E}(\widehat H) \right)$$
 and
  $$\tau_2 (H) = \nu \left( {\cal E}(\widehat H) \cap {\cal E}(\widehat H_{sim}) \right) / \nu \left( {\cal E}(\widehat H_{sim}) \right).$$
Results for all combination of parameters $(p(x),r)$ are reported in Table~\ref{tab:coverrate}.
The coverage rates are also computed according to the distance and plotted in Figure~\ref{fig:inter_009} ($\tau_1$ in solid line and $\tau_2$ in dashed line).
These results show that for any configuration and interaction range the approximation procedure of the conditional intensity of parent points in $W \cup \dW$ given the offspring points in $W$ is conservative.

\begin{table}[h]
    \centering
    \scalebox{0.95}{
    \begin{tabular}{|c|rrr|}
    \cline{2-4}
    \hline
   $r$ & $0.05$ & $0.09$ & $0.13$ \\
         \hline
         \rowcolor{Gray}
          \multicolumn{4}{c}{$p_1(x)$} \\
         \hline
$\tau_1 (H) $ & 88.04 &  93.48 &  98.73 \\
$\tau_2 (H) $ & 53.66 & 64.82 & 69.81  \\
\hline
$\tau_1 (E) $ &  97.30 &   97.40 &  99.91   \\
$\tau_2 (E) $ &  67.75 & 71.43 & 74.68  \\
\hline
$\tau_1 ( B^{inner}) $ & 100.00 &  99.86 & 100.00  \\
$\tau_2 (B^{inner}) $ & 87.75 & 78.21 & 81.69 \\
\hline
$\tau_1 (B^{outer}) $ & 100.00 & 100.00 & 100.00 \\
$\tau_2 (B^{outer}) $ & 81.63 & 72.37 & 74.30 \\
\hline
         \rowcolor{Gray}
  \multicolumn{4}{c}{$p_2(x)$} \\
\hline
$\tau_1 (H) $ & 90.07 &  96.83 &  99.37   \\
$\tau_2 (H) $ & 58.31 & 68.22 & 77.97 \\
\hline
$\tau_1 (E) $ &  97.72 & 96.24 & 90.65  \\
$\tau_2 (E) $ &  70.11 & 69.63 & 79.14  \\
\hline
$\tau_1 (B^{inner}) $ & 100.00 & 100.00 & 100.00 \\
$\tau_2 (B^{inner}) $ & 84.36 & 88.02 & 85.92 \\
\hline
$\tau_1 (B^{outer}) $ & 100.00 &  100.00 &  98.57  \\
$\tau_2 (B^{outer}) $ & 78.94 & 70.43 & 69.06 \\
\hline
    \end{tabular}}
    \caption{Coverage rate between the envelopes of the interaction statistics computed from the observed parent points and from the simulated parent points.}
    \label{tab:coverrate}
\end{table}

\section{Conclusion}

In order to quantify discrepancies between true conditional intensities and estimated ones (as in \cite{gabriel2021}), we have to both know the conditional intensity and to get fast computations to browse the space of conditioning realisations.
Because existing methods are computationally intensive, not allowing many simulations, instead of simulating the conditional intensity we proposed to consider an approximating process, whose deviation to the true process can be controlled.
We thus propose this approximation if one needs a fast procedure, even if conservative.

\bibliographystyle{apalike}
\bibliography{RefNote}

\end{document}